\documentclass[a4paper,10pt]{amsart}
\usepackage[utf8x]{inputenc}

\usepackage{cmap}
\usepackage{amsmath,amsfonts}
\usepackage[pdftex,unicode]{hyperref}
\usepackage{graphicx}

\newtheorem{thm}{Theorem}[section]
\newtheorem{lm}[thm]{Lemma}
\newtheorem{cor}[thm]{Corollary}
\newtheorem{defn}[thm]{Definition}

\newcommand{\Z}{Z\mathbb C[G]}

\newcommand{\tr}{\mathrm{tr}}
\title{Covers counting via Feynman calculus}
\author{Maksim Karev}
\thanks{Research is supported in part by the SNSF grant nr. 125070, 140666
and the TROPGEO project of the European Research Council.}
\address{Universit\'e de Gen\`eve, Section de Math\'ematiques, villa Battelle,
7, route de
Drize, 1227 Carouge, Switzerland}
\email{maksim.karev@unige.ch}
\begin{document}

\begin{abstract}
Let $G$ be a finite group. In this paper we present a tool for counting the
number of principle $G$-bundles over a
surface. As
an application, we express (non-standard) generating functions
for double Hurwitz numbers as integrals over  commutative Frobenius algebras,
associated with  symmetric groups.
\end{abstract}

\maketitle
\section{Introduction}

This note is devoted to a partial case of the following big problem.

\emph{Big problem:} Let $M=\{\mu_1,\ldots,\mu_k\}$ be a finite ordered
collection of 
conjugacy classes of a finite group $G$. Count the weighted number 
of principle $G$-bundles over a closed oriented surface $\Sigma$ of genus
$g$ with $k$ marked points $q_1,\ldots,q_k$, such that the
holonomy around $q_i$
is in the class $\mu_i$ for $i=1,\ldots,k$.

Here the ``weighted'' means that we count each bundle with a weight reciprocal
to
the number of its automorphisms.

It is the problem of the study of a certain 2-dimensional topological
quantum field theory (see, for example~\cite{A}). The corresponding numbers are
correlation functions in this theory.

An efficient way to deal with it is to use its connection with
the combinatorics and representation theory of the  group $G$. Namely, if
$\#N(M)$ is the number of
homomorphisms of the fundamental
group of $\Sigma/\{q_1,\ldots,q_k\}$ to $G$, such that the image of an
element corresponding to the complete revolution around $q_i$ is contained in
the conjugacy class corresponding to $\mu_i$, then the numbers of bundles is
the ratio of $\#N(M)$ and the order of $G$. For the details see~\cite{LZ}.



We use this connection for counting the number of bundles of a special
kind. Namely, denote by
  $h_\tau (\mu,\nu;r)$ the weighted numbers of
principle $G$-bundles over a
sphere with a non-trivial holonomy around $r+2$ points $q_0,\ldots,q_{r+1}$,
such that the holonomy around point $q_0$ is in the conjugacy class $\mu$, the
holonomy around $q_{r+1}$ in in the conjugacy class $\nu$, and
the holonomy around  all the other points is in the fixed conjugacy class
$\tau$.

The main result of this note is is the following. We construct 
generating functions for the numbers $h_\tau(\mu,\nu;r)$ as integrals over the
center of the group algebra $\mathbb C[G]$. Namely, for an explicitly
constructable square matrix $A_\tau(\beta)$, we have
$$
\sum_{k=0}^\infty h_\tau (\mu,\nu;r) \beta^r =
\frac{1}{Z_\tau}\frac{\tr(f_{\mu^{-1}} f_\mu)
\tr (f_{\nu^{-1}}
f_\nu)}{\# G}\int_{Z\mathbb C[G]} z_\mu \bar z_\nu e^{-<A_\tau(\beta)z,z>} dm,
$$
where all the notation is explained in the text.  This presentation allows us
to express the generating functions for numbers $h_\tau(\mu,\nu;r)$ as
entries of the matrix $(A_\tau(\beta))^{-1}$, multiplied by certain constants.

As an corollary of this presentation, we find that the corresponding generating
functions are rational functions in $\beta$. We also derive a non-linear
differental equation for this functions.

The main example for our considerations are (disconnected) double Hurwitz
numbers. It is the case when the group $G$ is a symmetric group $S_d$, and the
distinguished class $\tau$ is the class of a transposition. Principle
$S_d$-bundles are just ramified covers of degree $d$. 
These numbers have a very rich structure related with various fields of
mathematics.

In the case when one of the holonomies (say, $\nu$) is in the
conjugacy  class of  the unit, the celebrated ELSV-formula (see, for example the
original paper~\cite{ELSV})
relates the numbers of connected covers with intersection numbers on moduli
space of 
complex curves. 

The geometrical meaning of  double Hurwitz numbers  is not yet completely
understood. In the paper~\cite{GJV},  the authors conjectured than if $\nu =
(d)$, there exists a moduli space $Pic_{g,n}$ such that this numbers can be
expressed in the terms of intersection numbers on some properly chosen
compactification. They have also found an expression for generating functions
for double Hurwitz numbers in the terms of Schur polynomials and found an
explicit formula for one-part double Hurwitz numbers (it is the case when
there is a complete branching over one of the special points). In the
paper~\cite{O} there was shown that a generating function for double Hurwitz
numbers is a $\tau$-function for Toda lattice hierarchy.

Later in~\cite{CJM}, it was found that the computation of double Hurwitz numbers
 can be carried out in the language of tropical geometry.
In~\cite{BBM} the authors give a method for computation of number of covers
with an arbitrary ramification. In a sence, our work is a continuation of the
research carried out in~\cite{BBM}. 

The paper~\cite{J} provides the following formula for a generating function for
the double Hurwitz numbers, generalizing the results of~\cite{GJV}. Let the
cyclic type of $\mu$ be $(\mu_i)_{i=1}^{l(\mu)}$, the cyclic type of $\nu$ be
$(\nu_j)_{j=1}^{l(\nu)}$. By $\varsigma (z)$ we denote the following function
$$\varsigma (z) = e^{z/2} - e^{-z/2}.$$
Than, for $\mu$ and $\nu$ satisfying certain condition (the pair ($\mu$, $\nu$)
is
contained in a \emph{chamber} of $R_{l(\mu),l(\nu)}$, for the details
see~\cite{J}) we have
$$\sum_{r=0}^\infty h_\tau(\mu,\nu;r)\frac {z^r}{r!} = \frac 1{(\# Aut\ \mu)
\prod
\mu_i}\frac 1{(\# Aut\ \nu) \prod \nu_j} \frac 1{\varsigma(dz)}
\sum_{k=1}^{t(\mu,\nu)}\prod_{m=1}^{l(\mu)+l(\nu) +1} \varsigma
(zQ^{\mu,\nu}_{k,l}),$$
where $t(\mu,\nu)$ is a finite number, and $Q^{\mu,\nu}_{k,l}$ are certain
quadratic polynomials in $\mu_i$ and $\nu_j$.

This result is obtained using the the infinite wedge space formalism.

In the language of the infinite wedge space, the application of our technique
to the problem of double Hurwitz numbers computation is rather trivial. Namely,
we study the vacuum expectations of operators of the form
 $$\big\langle \prod \alpha_{\mu_i} (1 - \beta \mathcal F_2)^{-1} \prod
\alpha_{-\nu_j} \big\rangle,$$
which is nothing more but a study of entries of the block-diagonal operator
$(1-\beta \mathcal F_2)^{-1}$ in the basis $\prod \alpha_{-\nu_j}| 0 \rangle$.

Our technique allows to produce generating
functions for double Hurwitz numbers even in the case when a pair $(\mu,\nu)$
belongs to a \emph{resonance arrangement} in the terminology of~\cite{J}. Also
it works in the case of an arbitrary finite group.





The paper is organized as follows.  A
presentation of a fundamental group of a surface is constructed in the section
2. Section 3 presents a
reformulation of the enumerational problem in the algebraic language. Section
4 is a short overview of Feynman calculus, and section 5 is devoted to the
application of this technique to our problem.
Section 6 contains the calculation of (disconnected) double Hurwitz numbers for
the degrees $d=2,3,4$. We
also give a comparison of our generating function with the
generating function of~\cite{J} in one particular case.

The author expresses his sincere gratitude to J.~Rau and G.~Borot, who have
carefully read a preliminary version of the paper
paper, pointed out some inaccuracies and helped to make it more readable. I
also thank
G.~Mikhalkin, S.~Duzhin, P.~ Mnëv,
P.~Putrov, and N.~Kalinin for fruitful discussions concerning the subject.

\section{Group presentation corresponding to a pair-of-pants decomposition of a
surface}\label{pres}

We begin with a special presentation of the fundamental group of a punctured
surface.

\begin{defn}
 Let $g$ and $n$ be two non-negative integer numbers such that $2-2g-n < 0$
Denote by $H_{g,n}$ the fundamental group of a connected oriented surface
of genus $g$ with $n$ punctures.
\end{defn}

 Let us define some notions we are
going to use for construction of a presentation of $H_{g,n}$.

\begin{defn}
 By a \emph{1-3-valent graph} $\Gamma$ we mean a finite graph having only
vertices of valence 1 and 3. The first Betti number of $\Gamma$ (considered as a
one dimensional cell complex) is referred to as the \emph{genus}.
 The 1-valent vertices of a graph are
referred to as \emph{ends}; the edges, incident to ends, are called
\emph{leaves}; all the other edges are called \emph{inner edges}.  The set of
all vertices of $\Gamma$ which are not  ends is called the \emph{set of inner
vertices} and is denoted by $V^0(\Gamma)$. The set of all edges of $\Gamma$ is
denoted by $E(\Gamma)$; the set of inner edges of $\Gamma$ is denoted by
$E^0(\Gamma)$. The set of the edges of $\Gamma$ adjacent to $v\in V^0(\Gamma)$ is denoted by $E_v$.
\end{defn}

From now on, all the graphs we are dealing with, are supposed to be 1-3-valent,
if the contrary is not stated explicitly.

Another notion we are going to use is the following.

\begin{defn}
 Let $\Gamma$ be a connected graph. A \emph{maximal tree} $T$ in $\Gamma$ is a
connected subgraph of $\Gamma$ of genus $0$, such that the set of vertices of
$T$ coincides with the set of vertices of $\Gamma$.
\end{defn}

Each connected graph contains a maximal tree.

\begin{defn}
An  \emph{enhanced graph} is a connected graph enhanced with the following
additional data:\begin{itemize}
\item A choice of an orientation on all the edges of $\Gamma$; 
\item A choice of a cyclic order on a set of half-edges, adjacent to every 3-valent vertex of $\Gamma$;
\item A choice of a maximal tree $T$ in $\Gamma$;
\item A choice of a basepoint $p$ in $T$.
\end{itemize}
For a given edge $e\in E(\Gamma)$, and a given orientation on $e$, a vertex $v$
adjacent to $e$ is referred to as a \emph{source} of $e$, if the orientation of
$e$ is directed outwards with respect to $v$. If the orientation of $e$ is
directed inwards with respect to $v$, the vertex $v$ is referred to as a
\emph{sink} of $e$.
\end{defn}

We use the same notation for an enhanced graph and its underlying graph if it does not lead to an ambiguity.

The
starting point of the following construction is an enhanced  graph
$\Gamma$ with $n$ ends and the genus equal to $g$. Our construction of a presentation of $H_{g,n}$ consists of the following steps.

\begin{enumerate}

\item
Consider $\Gamma$ as a one-dimensional cell complex.  Subdivide each inner edge
of $\Gamma$ by a 2-valent vertex. Glue an additional 1-cell, in such a way that
both boundary components get glued to the same point, to each of the newly added
2-valent vertices, and to each of the ends of $\Gamma$. All the newly added
1-cells are referred to as \emph{circles}. Choose an orientation on all the
circles in an arbitrary way.
The obtained cell complex is homotopy equivalent to a bouquet of
$|E(\Gamma)|+g$ circles. Its fundamental group is a free group of the
corresponding rank. We present it in the following way.
\begin{itemize}
 \item \emph{Generators corresponding to edges in $T$} are presented by based
loops that start in $p$, go
along edges of $T$ to the  point of the attachment of  the circle corresponding
to
the edge we are interested in, make a complete  revolution along it in the
positive direction, and return back to $p$ along edges of $M$.
 \item \emph{Generators corresponding to edges that are not in $T$} are
presented by based loops with the base $p$.
According to the chosen orientation, each edge that does not belong to $T$ has
a beginning and an end. Start in $p$, go along edges of $T$ to the
beginning of the edge we are interested in, go along this edge to the point of
attachment of the circle, make a full revolution in the positive direction,
and return back the same way.
\end{itemize}

If $e$ is an edge of $\Gamma$, the described  above generator corresponding to
$e$ is denoted $p_e$.

\begin{itemize}
 \item \emph{Generators corresponding to  nontrivial cycles in $\Gamma$} are
presented by based loops with base $p$. The set of such generators is in the
bijection with the set of edges of $\Gamma$ that does not belong to $T$.  Start
in $p$, go along edges of $M$ to the
beginning of the edge which is not in $T$, continue
along it to the end, and return back to $p$ along edges of $T$.
\end{itemize}
If $e$ is an edge of $\Gamma$ which is not contained in  the maximal tree, the
above described generator is denoted by $g_e$.

\item
For each of the 3-valent vertices of $\Gamma$ prepare a 2-cell with an oriented
boundary.  Attach them to the cell complex in such
a manner that the boundary goes in the following way. We begin in a 3-valent
vertex, go along one of the edges attached to it to the point of the attachment
of the circle, and go along it using the following convention: if we came there
along the orientation of the edge, we go along the circle in the positive
direction. Otherwise, we go along it in the negative direction. Than we return
back to the 3-valent vertex we have started with and continue along a next
half-edge (according to the chosen cyclic order on the set of half-edges) adjacent to it. Repeating the
procedure for every adjacent edge, we obtain the way each of the  2-cells is
glued to the cell complex. For an example of such a gluing, see the figure~\ref{glue}.
\begin{figure}[ht]
\includegraphics[viewport = 1 0 170 215,scale=0.5]{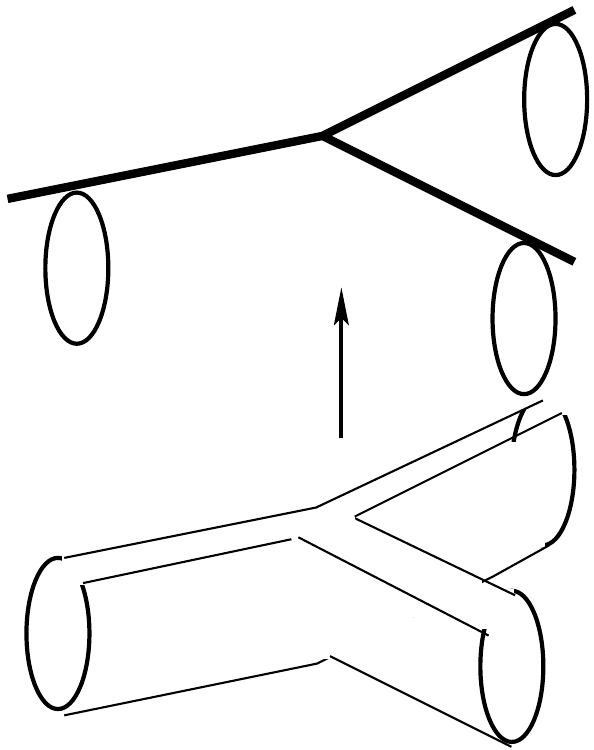}
\caption{An example of gluing a 2-cell to a graph with attached circles. A vicinity of a 3-valent vertex of a graph is shown in the top part of the picture.}\label{glue}
\end{figure}

The obtained cell
complex is homotopy equivalent to a connected oriented surface of genus
$g$ with $n$ punctures.  Each of the attached 2-cell provides a relation
in its fundamental group. Note, that for every edge $e$ not containing in the
maximal tree, the generator $g_e$ appears in only one of the obtained relations.
Van Kampen theorem implies that these relations are
defining relations in $H_{g,n}$. 
\end{enumerate}

Euler characteristics of a genus $g$ connected graph is $1-g$. It allows us to
compute explicitly the number of generators and relations in the presentation in
terms of the genus and the number of branching points. Therefore, we have just
proved the following Lemma.

\begin{lm}
The group $H_{g,n}$ admits a presentation with $4g - 3 + 2n $
generators and $2g - 2 + n$ relations. 
\end{lm}

Each 2-cell used in the construction corresponds to a pair of pants. So the obtained
presentation corresponds to a pair-of-pants decomposition of a surface.

\begin{figure}[hb]\label{presentation}
\includegraphics[viewport = 0 0 133 144,scale=0.7]{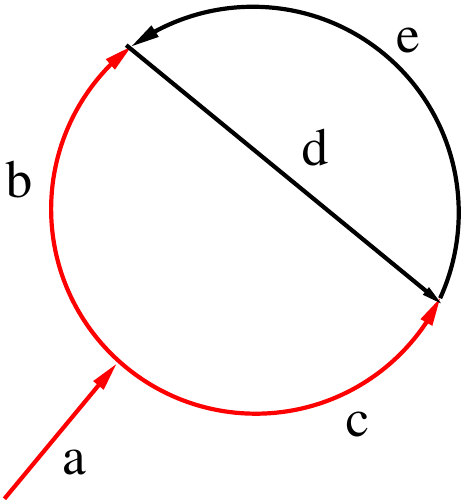}
\caption{An enhanced graph. The maximal tree is shown in red. The orientation of the edge is shown by an arrowhead. The cyclic ordering of half-edges obey the blackboard convention. The basepoint is the only 3-valent vertex of the maximal tree.}
\end{figure}

\textit{Example:} Let us consider the enhanced graph on the figure~\ref{presentation} with the blackboard cyclic ordering of half-edges in its vertices. The maximal tree is chosen
as it is depicted. The basepoint is the only 3-valent vertex in the maximal tree.

The corresponding group presentation is the following:
$$H_{2,1} = \langle p_a,p_b,p_c,p_d,p_e,g_d,g_e \,|\, p_a^{-1} p_b p_c\ =\ p_b^{-1}g_e^{-1}p_eg_ep_d\ = \
p_c^{-1}g_d^{-1}p_dg_dp_e\ =\ 1 \rangle.$$

\section{Problem reformulation in the terms of group algebra}

The problem of principle $G$-bundles
enumeration is
equivalent to the problem of enumeration of morphisms from the fundamental group
of a punctured surface to a finite group $G$ (for the details see, for
example~\cite{LZ}).

The number of possible morphism between two groups does not
depend on their presentations. So, considering the presentation constructed in
the previous section, we split the main problem of this note in two parts:
\begin{itemize}
 \item Fix the conjugacy classes of \emph{all} the images of the generators
of $H_{g,n}$ corresponding to the edges of the chosen graph $\Gamma$, and count
the number of such morphisms.
 \item Take a sum over all possible choices of conjugacy classes of images of
the inner edges. 
\end{itemize}

A similar method was discussed in~\cite{BBM}, where the authors consider the
problem from the tropical point of view. In a sense, from one hand, the main
result of this
paper is the generalization of the result of~\cite{BBM} on the case of an
arbitrary finite group. From the other hand, it is nothing else but a
well-known argument, which can be found, for example, in~\cite{D}, but in a
slightly non-standard form. That is why we present a complete proof.

To deal with the first part we need some algebraic notions. 

\begin{defn}
Let $G$ be a finite group. By a \emph{complex group algebra} of $G$ we mean the
algebra of complex-valued functions on $G$ with the convolution as a
multiplication. Its center is denoted by $Z\mathbb C[G]$. The set of conjugacy
classes of $G$ is denoted by $\Lambda_G$.
\end{defn}

The center $Z\mathbb{C}[G]$ of complex group algebra
has a natural basis indexed by $\Lambda_G$ (see, for example~\cite{LZ}, appendix
A).  Namely, if $\mu$ is
an arbitrary conjugacy class, the corresponding basis element is its characteristic function (here we regard a conjugacy class as a subset of $G$).
For a fixed group $G$ denote the basis element of $Z\mathbb{C}[G]$ corresponding
to the conjugacy class $\mu$ by $f_\mu$. This basis is referred to as the
\emph{standard basis}. The trace function $\tr_G$,
endowing $Z\mathbb{C}[G]$ with the Frobenius structure, is defined as follows.
\begin{defn} Let $G$ be a finite group. The \emph{trace function} is a linear functional $\tr_G \colon Z\mathbb C[G] \to \mathbb C$  defined on the elements of the standard basis as follows.
$$\tr_G(f_\mu)=\left\{\begin{array}{l} 1,\mbox{ if } \mu = (1),\\
										
						0, \mbox{
otherwise,}\end{array}\right.$$
where $(1)$ denotes the conjugacy class of the unit element in $G$.
\end{defn}

We omit the subscript $G$ in the notation for the trace function  when it does
not cause an ambiguity.

\begin{defn}
 Let $\Gamma$ be a graph and $G$ is a finite  group. By a
\emph{coloring of $\Gamma$} we mean a map $c\colon E(\Gamma) \to \Lambda_G$. If
$\Gamma$ is an enhanced graph, for a vertex $v\in V^0(\Gamma)$, edge $e\in E_v$
and coloring $c$, by $\bar c_v(e)$ we mean $c(e)$ if $e$ is oriented in such a
way that $v$ is a source for it, and the reciprocal conjugacy class otherwise.
\end{defn}

The procedure of the Section~\ref{pres} allows to construct a presentation  of a
fundamental group of a three times punctured sphere with 3 generators and 1
relation: use a enhanced graph $D$ with one 3-valent vertex, three 1-valent
ones. Fix such a presentation, and denote by $p_e$ the generator corresponding
to the edge $e\in E(D)$. Let $G$ be an arbitrary finite group, and $c$ is a
coloring of $D$. Denote by $N_D (c)$ the set of morphisms $H_{0,3}\to G$, such
that the image of the generator $p_e$ is contained in the conjugacy class $c(e)$
for all $e\in E(D)$. 

A tautological corollary of the definition of trace function reads:

\begin{lm}\label{triple}

Let $G$ be a finite group, and $c$ is a coloring of  an enhanced graph $D$
described above. Than the set of morphisms from $H_{0,3}$ to $G$, such that the
image of the generator $p_e$ corresponding to the edge $e\in E(D)$ is contained
in $c(e)$ has the cardinality
$$\# N_D (c) = \tr (\prod_{e\in E(D)}f_{\bar c_v(e)} ),$$
where $v$ denotes the only vertex of $D$.
\end{lm}

Now we are ready to prove the following Lemma.

\begin{lm}
Let $G$ be a finite group. Let $\Gamma$ be an enhanced graph  with coloring $c$.
Then the set $ N_\Gamma(c)$ of morphisms
$H_{g,n} \to G$, such that the image of a generator $p_e$ corresponding  to the
edge $e \in  E(\Gamma)$
belongs to the conjugacy class $c(e)$, has the cardinality
$$\# N_\Gamma (c) = |G|^g\big( \prod_{v \in V^0} \tr  (\prod_{e'\in E_v} f_{\bar
c_v(e')})\big) \prod_{e \in E^0} \frac 1{\tr (f_{c(e)} f_{c^{-1}(e)})},
$$
where $c^{-1}(e)$ denotes the conjugacy class reciprocal to
$c(e)$ for every edge $e\in E(\Gamma)$.
\end{lm}

\textit{Proof:} First of all, note, that for every conjugacy class $\mu\in \Lambda_G$, $\tr
(f_\mu f_{\mu^-1}) $ equals the number of group elements belonging to $\mu$.

Consider the enhanced graph  $\Gamma'$ which is  obtained from $\Gamma$ by
cutting all the edges not belonging to the maximal tree and attaching 1-valent
vertices to the remaining half-edges. The newly obtained leaves inherit the
orientation of edges and the coloring from cut edges. The cyclic orientation in
of half-edges adjacent to vertices remains the same. Denote the inherited
coloring by $c'$. The enhanced graph $\Gamma'$ allows to construct a
presentation of $H_{0,g+n}$.

The maximal tree for $\Gamma'$ coincides with $\Gamma'$.

The number $\# N_{\Gamma'}(c')$ is computed inductively.  Pick an arbitrary
vertex $v$ and make a random choice of images of the generators $p_e$
corresponding to the edges $e \in E_v$. According to Lemma~\ref{triple} it can
be made in $ \tr (\prod_{e\in E_v} f_{\bar c(e)})$ different ways.

Let  $v'$ be the second vertex adjacent to $e$. As the choice of  the image for
the generator $p_e$ is already made, the choice for the images for the
generators corresponding to the remaining edges adjacent to $v'$ can be
performed in $$\frac {\tr (\prod_{e'\in E_v'} f_{\bar
c(e')})}{\tr(f_{c(e)}f_{c^{-1}(e)})}$$ ways.

Any two vertices in $\Gamma'$ are connected with at most one edge.  It implies
that
$$\# N_{\Gamma'} (c') = \prod_{v \in V^0(\Gamma)}  \tr (\prod_{e'\in E_v}
f_{\bar c_v(e')}) \prod_{e \in E^0(\Gamma')} \frac 1{\tr (f_{c(e)}
f_{c^{-1}(e)})}.$$

 There is a natural map $b\colon E(\Gamma') \to E(\Gamma)$, which sends  an edge
of $\Gamma'$ to the corresponding edge of $\Gamma$ before cut. 
 The map $b$ gives rise to a morphism $\phi\colon H_{0,g+n}\to H_{g,n}$. On the level of generators, this morphism is described in  the following way:
\begin{itemize}
 \item If $e'$ is an inner edge of $\Gamma'$, or $e'$ is a leaf that is shared 
by both $\Gamma$ and $\Gamma'$, the generator $p_{e'}$ maps to the generator
corresponding to the edge $b(e')\in E(\Gamma)$.
 \item If $e'$ is a leaf of $\Gamma'$ came from a cut edge of $\Gamma$, and  its
end is a sink of $e'$, the generator $p_{e'}$ maps to the generator
corresponding to the edge $b(e')\in E(\Gamma)$.
 \item If $e'$ is a leaf of $\Gamma'$ came from a cut edge of $\Gamma$, and  its
end is a source of $e'$, the generator $p_{e'}$ maps to the $g^{-1}_{b(e')}
p_{b(e')} g_{b(e')}$.
\end{itemize}

The check that $\phi$ is indeed a group morphism is straightforward.

The morphism $\phi$ induces a natural map $\phi^*\colon N_\Gamma (c) \to N_{\Gamma'} (c')$.

For an edge $e\in E(\Gamma)$ that is not contained in the  maximal tree denote
by $e^{+}$ its preimage such that its end is the sink of $e^{+}$. The other
preimage is denoted by $e^-$.

Let $f'\in N_{\Gamma'} (c')$. We can construct out of it a morphism $f\in N_\Gamma (c)$ in the following way. 
If $e\in \Gamma$ is contained in the maximal tree, the image of $f(p_e)$
coincides  with the image $f'(p_{b^{-1}(e)})$. Otherwise, the image $f(p_e)$
coincides with $f'(p_{e^+})$.

 The only thing that is left is to fix images of generators $g_e$ corresponding
to  edges of $\Gamma$ which are not contained in the maximal tree. We fix them
in such a way, that  $f'(p_{e^-}) = f(g_e^{-1})f'(p_{e^+}) f(g_e)$. The check
that $f$ is indeed a group morphism is straightforward.

Moreover, it is clear that $\phi^*(f) = f'$. As for every $e$ not belonging to
the maximal tree, the choice of the image of $g_e$ can be performed in
$$\frac{|G|}{\tr(f_{c(e)} f_{c^{-1}(e)})}$$ ways, and the choices  for
different edges can be made completely independently, as  each $g_e$ appears in
only one relation, we arrive to the assertion
of the Lemma.
\nopagebreak

\hfill\textit{Q.E.D.}

Recall, that in fact we are interested in the number of morphisms $H_{g,n} \to
G$ such that only the conjugacy classes of images of the generators
corresponding to the leaves are fixed.

\begin{defn} Let $G$ be a finite group, $\Gamma$ be a graph. A
\emph{boundary condition} is a  map $M \colon E(\Gamma)\backslash E^0(\Gamma)
\to
\Lambda_G$. The set of all possible colorings $c\colon E(\Gamma) \to \Lambda_G$
such that $M = c|_{E(\Gamma)\backslash E^0(\Gamma)}$ is denoted by $C(M)$.
Abusing 
the notation denote $\cup_{C(M)} N_\Gamma (c) = N_\Gamma (M)$.
\end{defn}

\begin{thm}\label{main}
Let $G$ be a finite group, $\Gamma$ be an enhanced graph, and $M$ be a
boundary condition. Than the number $\# N_\Gamma (M) = \sum_{c\in C(M)} \#
N_\Gamma (c)$ is invariant under the following enhanced graphs transformations:
\begin{itemize}
 \item Other choice of the maximal tree in $\Gamma$;
 \item Change of orientation of an edge $e\in E^0(\Gamma)$;
 \item Change of cyclic order of half-edges adjacent to a vertex $v\in V^0(\Gamma)$;
 \item Local transformation of $\Gamma$ as it is depicted on the
figure~\ref{ihx}.
\end{itemize}
\begin{figure}[ht]
\includegraphics[viewport = 0 0 385 184,scale=0.5]{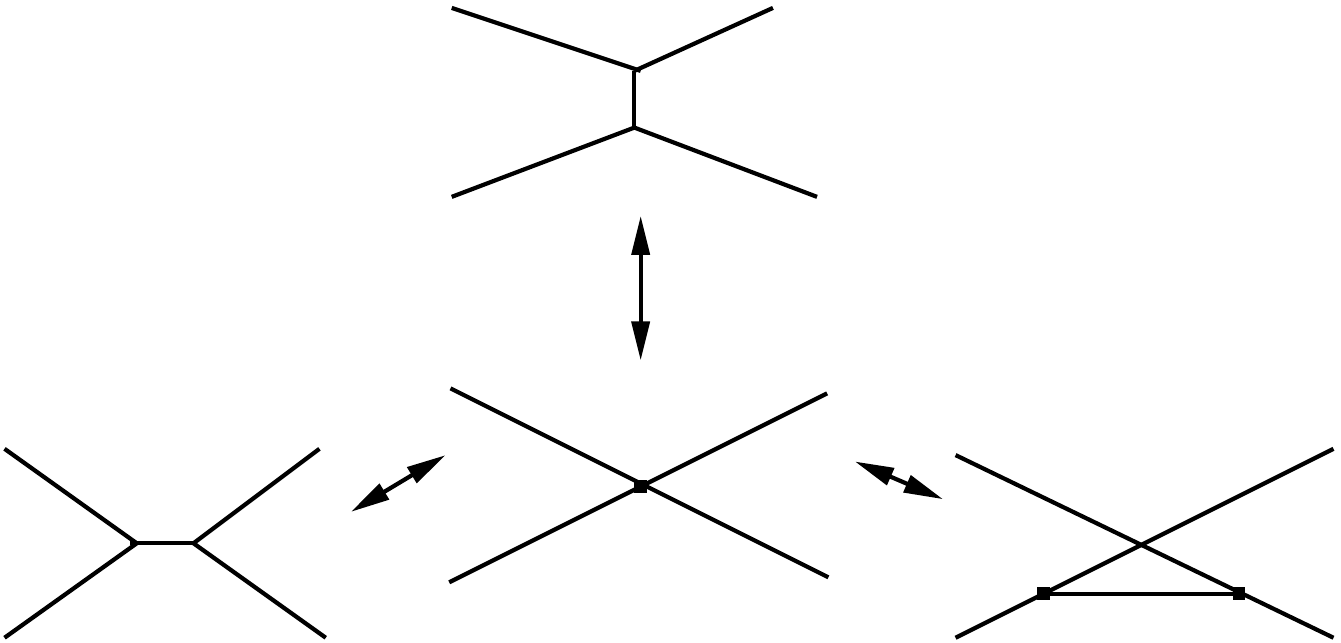}
\caption{The local graph transformations. The figure shows the possible resolutions of a 4-valent vertex.}\label{ihx}
\end{figure}
\end{thm}
 
\emph{Proof:}
A different choice of the maximal tree in $\Gamma$ respect every single summand
in the sum $\sum_{c\in C(M)} \# N_\Gamma (c)$.

Suppose, an edge  $e\in E^0(\Gamma)$ has a coloring $c(e)$. A change of the
orientation of an inner edge $e$, and switching its coloring to the reciprocal
$c^{-1}(e)$ also preserves every single term
in the sum under consideration. 

The independence  on the choice of the cyclic ordering of half-edges adjacent to
a vertex $v\in V^0(\Gamma)$ follows from the commutativity of 
$Z\mathbb C[G]$.

The  invariance under the local transformation of $\Gamma$ follows from the fact
that
the number of possible morphisms from one group to another (with some fixed
conditions) does not depend on
the presentation on the group. But we would like to give an alternative argument.

The commutative algebra $Z\mathbb C[G]$ can be endowed with a Hermitian product.
Namely, consider a
semilinear conjugation acting on a basis element $f_\mu$ for $\mu \in
\Lambda_G$ in the following
way:
$$\bar {\cdot} \colon f_\mu \mapsto f_{\mu^{-1}}.$$
The product is defined by formula $(a,b) = \tr(a\bar b)$ for $a,b\in Z\mathbb
C[G]$. The basis vectors of $Z\mathbb C[G]$ form an orthogonal basis with
respect to this product.

Consider a trace of a product of four basis vectors of $Z\mathbb C[G]$
corresponding to conjugacy classes $\mu_1,\mu_2,\mu_3,\mu_4\in \Lambda_G$. Due
to
general properties of Hilbert spaces, and commutativity of $Z\mathbb
C[G]$, for any permutation $\sigma \in S_4$, this trace can be expanded:
$$\tr(f_{\mu_1}f_{\mu_2}f_{\mu_3} f_{\mu_4}) = \sum_{\nu \in \Lambda_G}
\frac{\tr( f_{\mu_{\sigma(1)}}f_{\mu_{\sigma(2)}}f_\nu) \tr(
f_{\mu_{\sigma(3)}}f_{\mu_{\sigma(4)}}f_{\nu^{-1}})}{\tr (f_\nu
f_{\nu^{-1}})}.$$
This formula implies the desired invariance.
\nopagebreak

 \hfill\textit{Q.E.D.}

As the space of all connected 1-3-valent graphs is connected with respect to
above mentioned transformations (see, for example~\cite{H}), the only data that
matters for computation of the number $\# N_\Gamma (M)$ is a genus of an
enhanced graph $\Gamma$, the orientation of the leaves of $\Gamma$, and boundary
condition $\mu$. If we are dealing with a group, where every conjugacy class is
a self-reciprocal, like a symmetric group $S_n$, the orientation of the
leaves is also not relevant at all.

\textit{Example:} Let $\Gamma(2,1)$ be a connected 1-3-valent graph of genus 2
with 1 leaf, and $(-1)$ be the conjugacy class of the non-unit element of
$S_2$. Provide $\Gamma$ with an arbitrary enhancement. We are interested in the
computation of the number $\# N_{\Gamma(2,1)}(M)$, where $M$ is a boundary
condition which assigns class $(-1)$ to the only leaf of $\Gamma(2,1)$. As it
was stated above, this number is independent on the choice of the enhancement of
$\Gamma$.

Let $c$ be a coloring of $\Gamma(2,1)$. The multiplication table of $Z\mathbb
C[S_2]$ imply that if there exists a vertex $v\in V^0(\Gamma(2,1))$, such that
there are exactly one or three edges adjacent to $v$ colored by $(-1)$, than
this coloring does not contribute to the sum.

By the contrary, if every vertex $v\in V^0(\Gamma(2,1))$ is adjacent to
exactly zero or two edges of coloring $(-1)$, this coloring gives a contribution
1.

This actually means, that $\# N_{\Gamma(2,1)}(M) = 0$.

From the other hand, for the same group $S_2$, and for a graph $\Gamma(2,2)$ of
genus 2 with 2 leaves, and the boundary condition $M'$ which assigns $(-1)$
to every leaf of the graph, the number $\# N_{\Gamma(2,2)}(M') $ is 16,
which coincides with the prediction of the Frobenius formula (see, for example~\cite{LZ}, appendix A).

\section{Complex Feynman Calculus}

This section can be considered as a short exposition of the theory of the
Feynman calculus.

The Feynman calculus is a powerful tool for enumerating graphs. We will apply
it to compute generating functions for number of principle $G$-bundles.
Briefly speaking, the machinery of Feynman calculus works as follows. Consider
we have a graph with edges colored by a finite number of colors. Feynman
calculus  provides us with a rule how to assign a number to every vertex and
every
edge of a such a graph. The weight of a graph is calculated as a product of the
above
mentioned numbers over all the vertices and edges divided by the order of the
automorphism group of the graph. Than we take a sum of weights over all possible
graphs. It turns out, that this sum can be interpreted as a result of a
computation of an integral. For the details on Feynman calculus in the real case
we refer the reader to~\cite{E} or~\cite{LZ}.

M.~Mulase and J.~Yu in~\cite{MY} studied the Feynman calculus over a von Neumann
algebra.  In fact, our considerations are a generalization of the methods of
Mulase and Yu on the case of graphs with 1-valent vertices. 

In a sense, our considerations can be treated as a theory of Feynman calculus
over a space of diagonal matrices.

Here we will be interested in the complex Feynman calculus. The reason for it
is that complex calculus works for any finite group $G$, while the real
calculus works only in the cases when every conjugacy class of $G$ is a
self-reciprocal class.

We will need some
preliminary considerations to proceed.

\begin{defn}
 Let $P$ be a finite-dimensional complex vector space endowed with a Hermitian
product $<\cdot,\cdot>$. . By
$\mathbb C^n$ we mean an $n$-dimensional complex vector space with a chosen
orthonormal basis $\{f_i\}_{i=1}^n$ with respect to the chosen Hermitian
product. The measure $dm$ is defined as $dm = \prod_{i=1}^n
d\mathrm{Re}z_i\ d\mathrm{Im}z_i$, where $z_i$ is the coordinate corresponding
to the basis element $f_i$.
\end{defn}

\begin{lm}\label{int}
 Let $A$ be an endomorphism of $\mathbb C^n$, expressed in a standard basis as
a positively definite symmetric real matrix, and $p$ be
an arbitrary vector in $\mathbb C^n$. The following equalities are satisfied:
\begin{enumerate}
 \item $\int_{\mathbb C^n} e^{ -<Az,z> } dm =
\frac{\pi^n}{\mathrm{det}\ A};$
 \item $\int_{\mathbb C^n} e^{ -<Az,z> + <p,z> + <z,p> } dm =
\frac{\pi^n}{\mathrm{det}\ A} e^{<A^{-1}p,p>}.$
\end{enumerate}
\end{lm}

\textit{Proof:}
The first equality can be deduced from the fact that every sesquilinear
positively definite form can be diagonalized by the action of $U(n)$.
So it is enough to verify this formula in one dimensional case, where it is
nothing else but a product of two Gaussian integrals.

The second equality follows from the first one by the variable change $z
\mapsto z + A^{-1}p.$ 
\nopagebreak

\hfill\textit{Q.E.D.}

The immediate consequence of the lemma are the equalities:
$$\int_{\mathbb C^n} z_i \bar z_j e^{ -<Az,z> } dm =
\frac{\pi^n}{\mathrm{det}\ A} A^{-1}_{ij}\ \mbox{for all}\ i,j,$$
$$\int_{\mathbb C^n} z_i  z_j e^{ -<Az,z> } dm =
0,$$
$$\int_{\mathbb C^n}\bar z_i \bar z_j e^{ -<Az,z> } dm =
0.$$
They can be obtained by the differentiation of the second equality of the Lemma by $p$ or $\bar p$.

\begin{defn}
 Let $k$ be a positive integer, and $K_1, K_2$ be two finite sets. By a
\emph{pairing} we mean a bijection $\sigma \colon K_1 \to K_2$. The set of all
pairings  is denoted by $\Pi(K_1,K_2)$.
\end{defn}

It is clear that $\#\Pi(K_1, K_2) = \delta_{\# K_1,\#K_2} (\#K_1)!$, where
$\delta_{\cdot,\cdot}$ is the Kronecker symbol.
 
The next theorem is a standard fact. Its proof in the real case can be found
in~\cite{E}. The proof for the complex case can be obtained by the same
considerations.

\begin{thm}[Wick]\label{wick} Let  $L=\{l_1,
\ldots, l_m\}$, $L'=\{l'_1,\ldots,l'_n\}$ be
two collections of complex-linear forms  on $\mathbb
C^n$, and $A$ be an endomorphism of $\mathbb C^n$ expressed in the standard
basis by a symmetric positively definite real matrix. Then
$$\int_{\mathbb C^n} \prod_{l \in L_1} l(z)  \prod_{l'\in L'} \bar l'(z)e^{
-<Az,z>
} dm =
\frac{\pi^n}{\mathrm{det}\ A}\sum_{\sigma \in \Pi(L,L')} \prod_{l\in L}
<A^{-1}l,\sigma(l)>,$$
and the integral converges absolutely.
\end{thm}

Let us again stress  the fact that the only way to obtain
non-zero is to have the same number of linear and semilinear forms under the
integral.

\section{Double Hurwitz numbers}\label{hur}

\begin{defn}
 Let $\Sigma$ be an oriented surface and $f,g\colon \Sigma' \to \Sigma$ be
two principle $G$-bundles with a non-trivial holonomy a finite number of
points. Than $f$
and $g$ are considered to be equivalent, if there exists an automorphism
$h\colon \Sigma'\to \Sigma'$ such that $f = g \circ h$. The \emph{automorphism
group} $Aut\ f$ of a principle $G$-bundle $f\colon \Sigma' \to \Sigma$ is a
group
of automorphisms $h\colon \Sigma' \to \Sigma'$ such that $f = f\circ h$. The
number $\frac 1{\# (Aut\ f)}$ is referred to as the \emph{weight} of the
cover.
\end{defn}

Consider the following enumerative problem.

\emph{Problem:} Fix $k+2$ points $q_0,\ldots,q_{r+1}$ on $\mathbb CP^1$ and a
finite group $G$. Choose two conjugacy classes $\mu$ and $\nu$ of $G$. Find the
weighted number principle $G$-bundles
$\mathbb {C}P^1$, with a non-trivial holonomy around $q_0,\ldots,q_{r+1}$, such
that the
holonomy around $q_0$ is in $\mu$, the holonomy around $q_{r+1}$ is in
$\nu$, and the holonomy around all the other points $(q_1,\ldots,q_n)$ belong
to a fixed conjugacy class $\tau$.

 Denote the corresponding number by
$h\tau(\mu,\nu;r)$.

As it was mentioned in the introduction, the  weighted number of principle
$G$-bundles with fixed conjugacy class of holonomies, can be computed as a
number of the morphisms from the fundamental
group of a $r+2$-times punctured sphere to the group $G$ divided
by the order of $G$. 

\begin{defn}
 Let $G$ be a finite group and $\{f_\mu\}_{\mu \in \Lambda_{G}}$ be the
standard basis of $\Z$. Endow $\Z$ with a Hermitian product $<\cdot, \cdot>_G$,
such that the basis $\{f_\mu\}_{\mu\in \Lambda_G}$ is orthonormal with respect
to it. Let $\tau \in \Lambda_{G}$ be a fixed conjugacy class. Denote by
$A_\tau(\beta)$ a linear endomorphism of $Z\mathbb
C[G]$ represented in the standard basis as the matrix with entries
$$(A_\tau)_{\mu,\nu} (\beta)= \tr(f_{\mu^{-1}} f_\nu (1-\beta f_\tau)) $$ for
$\mu,\nu \in
\Lambda_{G}$.
\end{defn}

Note, that for any $\tau$ there exists a neighborhood of the point $\beta=0$
such
that
the matrix $A_\tau(\beta)$ is non-degenerate. 

\begin{thm}
 Let $G$ be a finite group, and  $\mu$, $\nu$ and $\tau$ be three elements of
$\Lambda_{G}$. The generating function $h_{\mu,\nu}^{\tau}(\beta) = \sum_r
h_\tau(\mu,\nu;r) \beta^r$ is given by the following formula:
$$h^{\tau}_{\mu,\nu}(\beta) = \frac{1}{Z_\tau}\frac{\tr(f_{\mu^{-1}} f_\mu) \tr
(f_{\nu^{-1}}
f_\nu)}{\# G}\int_{Z\mathbb C[G]} z_\mu \bar z_\nu e^{-<A_\tau(\beta)z,z>}
dm,$$
where $Z_\tau$ denotes the value of the
integral
$$Z_\tau = \int_{Z\mathbb C[G]} e^{-<A_\tau(\beta) z,z>} dm.$$
\end{thm}

\textit{Proof:} 
Let $\Gamma_r$ be an enhanced graph of genus $0$ with $r+2$ leaves such that
every inner vertex of $\Gamma_r$ is adjacent to at least one leaf. It is clear,
that there is exactly two inner vertices $v_1$ and $v_2$ of $\Gamma_r$ adjacent
to two leaves. Denote one of the leaves adjacent to $v_1$ as $e_1$, and one of
the leaves adjacent to $v_2$ as $e_2$. Orient the leaves of $\Gamma_r$ such
that for each of them the corresponding end is a sink.

Denote by $M_r$ a boundary condition which assigns the class $\mu$ to $e_1$,
the class $\nu$ to $e_2$, and class of a permutation $\tau$ to every other leaf
of $\Gamma_r$. 

It is clear from our previous considerations that $h_\tau(\mu,\nu,r) = \frac
{\#N_{\Gamma_r}(M_r)}{\# G}$.

According to the Theorem~\ref{main}, the number of covers can
be computed by the formula:
$$h_\tau (\mu,\nu;r) = \frac 1{\# G} \sum_{c \in C(M_r)} \Big( \prod_{v \in
V_0(\Gamma_r)} 
\tr ( \prod_{e'\in E_v}
f_{\bar c_v(e')})\Big) \prod_{e \in E_0(\Gamma_r)} \frac 1{\tr (f_{c(e)}
f_{c^{-1}(e)})}.$$

 Let $k$
be a non-negative integer. Denote by $X_{r}$  the set of maps $$x\colon
\{1,\ldots,r\} \to (\Lambda_{S_d})^2.$$ The notation $x_i(k)$ for $i=1,2$
denotes the corresponding component of the map. By an automorphism of a map
$x\in X_r$ we mean an automorphism $h$ of the set  $\{1,\ldots,r\}$ such that $x
= x\circ h$. The group of automorphisms of $x\in X_r$ is denoted by $Aut\
x$. 

Let $D$ be a linear automorphism of $Z\mathbb C[G]$ represented in a standard
basis by a matrix $D_{\mu,\nu} = \tr (f_{\mu^{-1}} f_\nu)$.

Let us discuss the denominator $Z_\tau$ in the statement of the Theorem. Using
the
absolute convergence of the integral in a neighborhood of the point $\beta =
0$, we
expand the expression in the following manner:

$$Z_\tau = \sum_{r=0}^\infty \beta^r \sum_{x\in X_r}\frac{1}{\#Aut\ x}
\int_{Z\mathbb C[G]} \prod_{k=1}^{r}
\tr \big( f_{x_1(k)} f_{x_2(k)} f_\tau \big) z_{x_1(k)} \bar
z_{x_2(k)}e^{-<Dz,z>}dm.$$

Apply theorem~\ref{wick} for the computation of this expression. For $x\in
X_r$ let $\Pi^x$ denote a subset of $\Pi_r = \Pi (\{1,\ldots,r\},\{1,\ldots,
r\})$, such that for any $\sigma\in \Pi^x$ $x_1(k)= x_2^{-1}(\sigma(k))$. The
elements of $\Pi^x$ are called \emph{admissible pairings}. For any
non-negative integer $r$, we obtain

$$Z_\tau = \frac {\pi^{\# \Lambda_{G}}}{\mathrm{det}\ D} \sum_{r=0}^{\infty}
\beta^r \sum_{x\in X_r} \sum_{\sigma\in \Pi^x} \frac 1{\# Aut\ x}
\prod_{k=1}^{r} \frac {\tr \big( f_{x_1(k)} f_{x_2(k)} f_\tau
\big)} {  \tr(f_{x_1(k)}f_{x_2^{-1}(\sigma(k))})}.$$

This expression can be interpreted in the following way. Let  a
\emph{flower} be an elementary piece of a graph consisting of two vertices
connected by an edge and two oriented half-edges attached to one of the
vertices. For one of the half-edges the adjacent vertex is a sink (we call
this half-edge \emph{left}), for another one the adjacent vertex is a source
(this half-edge is \emph{right}). The edge connecting the
vertices of a flower is oriented, in such a way that its end is a sink. Each
half-edge of a flower is
colored by an element of $\Lambda_{G}$, the only edge of a flower is colored
by the fixed class $\tau$ of  $G$.

\begin{figure}[ht]
\includegraphics[viewport = 2 2 204 127,scale=0.4]{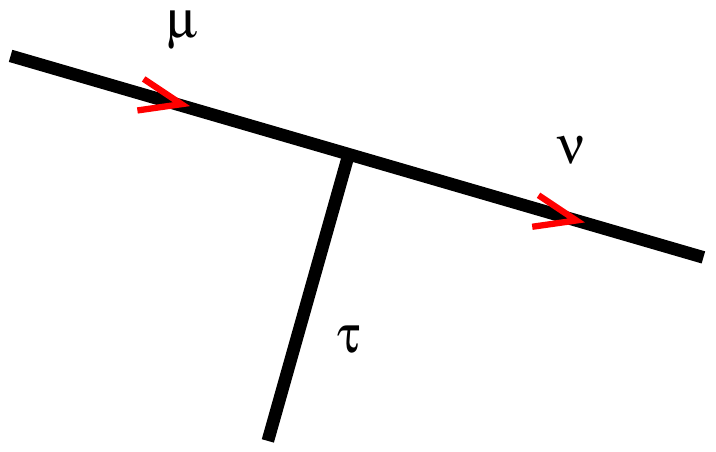}
\caption{A flower. The left half-edge is colored by $\mu\in \Lambda_{G}$, the
right half-edge is colored by $\nu\in \Lambda_{G}$.}\label{flower}
\end{figure}

Every pairing $\sigma \in \Pi_r$ corresponds to an element of the symmetric
group $S_r$. Namely, set the result of application of this element to $k\in
\{1,\ldots,r\}$ to be equal $\sigma(k)$. By the abuse of notation, denote a
pairing and the corresponding permutation by the same letter.

\begin{figure}
\includegraphics[viewport =1 1 245 242, scale=0.7]{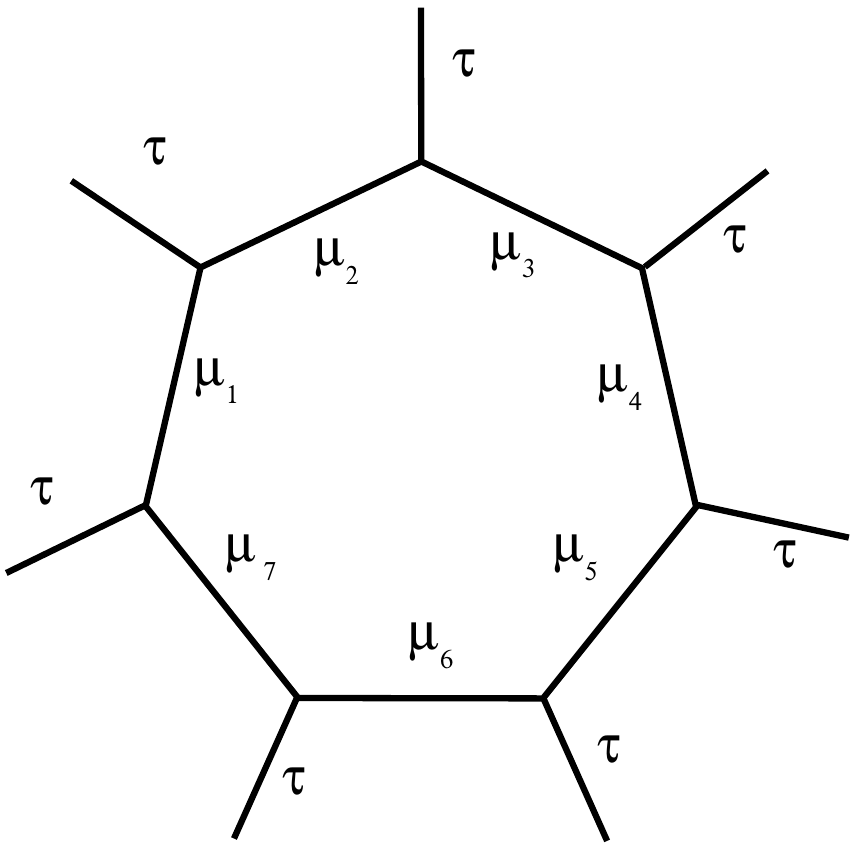}
\caption{An example of assembling a graph. The directions of the edges are not
shown.}\label{neck}
\end{figure}

Pick an arbitrary $x\in X_r$. For every $k\in \{1,\ldots,r\}$ prepare a flower
with the left half-edge colored by $x_1(k)$, and the right half-edge
colored by $x_2(k)$. Each admissible pairing provides a recipe to assembly a
colored 1-3-valent graph. Namely, the right half-edge of the $k$th flower gets
glued to the left half-edge of the $\sigma(k)$th flower. The resulting graph
can be disconnected, but each of its connected components has genus 1. The
number of connected components of the graph equals the number of cycles in the
permutation $\sigma$. For an example see the figure~\ref{neck}.

We claim, that Theorem~\ref{main} implies that $Z_\tau$ is a generating function
for the weighted number of principle $G$-bundles over a collection of disjoint
tori, such that every non-trivial holonomy around a point is in the class $\tau$
multiplied by $ {\pi^{\# \Lambda_{G}}}/{\mathrm{det}\ D}$. Note, that from the
other hand, $Z_\tau = {\pi^{\# \Lambda_{G}}}/{\mathrm{det}\ A_\tau(\beta)}$.

Now consider the integral
$$Z = \int_{Z\mathbb C[G]} z_\mu \bar z_\nu e^{-<A_\tau(\beta)z,z>} dm. $$
As it was mentioned during the study of $Z_\tau$, we can use the expansion:
$$ \sum_{r=0}^\infty \beta^r \sum_{x\in X_r}\frac{1}{\#Aut\ x}
\int_{Z\mathbb C[G]} z_\mu \bar z_\nu \prod_{k=1}^{r}
\tr \big( f_{x_1(k)} f_{x_2(k)} f_\tau \big) z_{x_1(k)} \bar
z_{x_2(k)}e^{-<Dz,z>}dm.$$

The only difference from $Z_\tau$ are two distinguished forms $z_\mu$ and $\bar
z_\nu$. Applying the Theorem~\ref{wick} again, we see, that the computation of
this integral differs from the computation of $Z_\tau$ only by the existence of
two additional elementary blocks. One is represented by a vertex with an
adjacent half-edge oriented in such a way that the vertex is a source, colored
by $\mu$. The other is represented by a vertex with an
adjacent half-edge oriented in such a way that the vertex is a sink, colored
by $\nu$. Our rules for the assembling graphs out of flowers and this two
additional details are the same. For each graph we pick only one copy of each
of these details.

Note, that both distinguished details always contribute to the same connected
component of a graph of genus 0. This connected component is exactly the graph
we discussed in the beginning of the proof.

The division by $Z_\tau$ allows to get rid of the contribution of all the
components of genus 1. The origin of the factor
${\tr(f_{\mu}^{-1} f_\mu) \tr (f_{\nu}^{-1}
f_\nu)}/{\# G} $ is obvious.
\nopagebreak

 \hfill\textit{Q.E.D.}

The direct application of the theorem~\ref{wick} lead to the following corollary:

\begin{cor}
$$h^{\tau}_{\mu,\nu}(\beta) = \frac{\tr(f_{\mu^{-1}} f_\mu) \tr (f_{\nu^{-1}}
f_\nu)}{\# G} (A^{-1}_\tau(\beta))_{\mu,\nu},$$
where $A^{-1}_\tau(\beta)$ is the matrix  inverse to $A_\tau(\beta)$.
\end{cor}

As the matrix $A_\tau(\beta)$ is linear in $\beta$, we obtain:

\begin{cor}
 For any finite group $G$ and any $\mu,\nu,\tau\in \Lambda_G$ the generating
function  $h^\tau_{\mu,\nu}(\beta)$ is a rational function in
$\beta$.
\end{cor}

Another simple theorem follows from our representation of these numbers in the
terms of graphs.

\begin{thm}
 The generating functions for numbers $h_\tau(\mu,\nu;r)$ satisfy the following
 equation
$$ h^{\tau}_{\mu,\nu}(\beta) + \beta \frac{\partial}{\partial \beta}
h^{\tau}_{\mu,\nu}(\beta)= (\# G) \sum_{\lambda \in \Lambda_{S_d}}
\frac{h^{\tau}_{\mu,\lambda}(\beta)h^{\tau}_{\lambda,\nu}(\beta)}{\tr
(f_{\lambda^{-1}} f_\lambda)}.
$$
\end{thm}

This equation corresponds to presentation of a genus 0 chain-like graph as a union of two graphs of the same type along their leaves.

\section{Double Hurwitz numbers}

Let $d$ be a positive integer, and $\tau_d\in \Lambda_{S_d}$ be the class of a
transposition.
In this section we present matrices $A_{\tau_d}(\beta)$ and
$(A_{\tau_d}(\beta))^{-1}$ for $d=2,3,4$.

For $d=2$ the basis elements are ordered in the following manner. The
first one corresponds to partition $(11)$, the second one corresponds
to $(2)$.
$$A_{\tau_2}(\beta) = \begin{pmatrix} 1&-\beta\\
                                                                                
                -\beta&1\end{pmatrix}\qquad (A_{\tau_2}(\beta))^{-1} =
\begin{pmatrix}
                                                                                
                \frac 1{1-\beta^2}&\frac{\beta}{1-\beta^2}\\
                                                                                
                \frac{\beta}{1-\beta^2}&\frac{1}{1-\beta^2}
                                                                                
                \end{pmatrix}.$$

For $d=3$ we use the ordering $(111),(12),(3)$.
$$A_{\tau_3}(\beta) = \begin{pmatrix}
1& -3\beta & 0\\
-3\beta & 3 & -6\beta\\
0 & -6\beta & 2
\end{pmatrix}\quad
(A_{\tau_3}(\beta))^{-1} = \begin{pmatrix}
\frac{6\beta^2-1}{9\beta^2-1} & -\frac{\beta}{9\beta^2-1} &
-\frac{3\beta^2}{9\beta^2-1}\\
-\frac{\beta}{9\beta^2-1}         & -\frac{1}{27\beta^2-3}     &
-\frac{\beta}{9\beta^2 -1}\\
-\frac{3\beta^2}{9\beta^2-1}   & -\frac{\beta}{9\beta^2-1} &
\frac{3\beta^2-1}{9\beta^2-1}
\end{pmatrix}
$$

For $d=4$ the ordering of basis elements is $(1111),(112),(22),(13),(4)$.
$$
A_{\tau_4}(\beta) = \begin{pmatrix}
1 & -6\beta & 0 & 0 & 0\\
-6\beta & 6 & -6\beta & -24\beta& 0\\
0 & -6\beta & 3 & 0 & -12\beta\\
0& -24\beta & 0 & 8 & -24\beta\\
0 & 0 & -12\beta & -24\beta & 6
\end{pmatrix}$$
And $(A_{\tau_4}(\beta))^{-1}$ is
$$\begin{pmatrix}
\frac{24\beta^4-34\beta^2+1}{144\beta^4-40\beta^2+1} &
-\frac{20\beta^3-\beta}{144\beta^4-40\beta^2+1} & \frac{24\beta^4 +
2\beta^2}{144\beta^4 -40\beta^2 +1} & -\frac{3\beta^2}{36\beta^2-1} &
\frac{16\beta^3}{144\beta^4-40\beta^2+1}\\
-\frac{20\beta^3-\beta}{144\beta^4-40\beta^2+1} &
-\frac{20\beta^2-1}{864\beta^4-240\beta^2+6} &
\frac{12\beta^3+\beta}{432\beta^4 -120\beta^2 +3} &
-\frac{\beta}{72\beta^2-2} & \frac{8\beta^2}{432\beta^4 -120\beta^2
+3}\\
\frac{24\beta^4 + 2\beta^2}{144\beta^4-40\beta^2+1} &
\frac{12\beta^3+\beta}{432\beta^4-120\beta^2+3} & \frac{72\beta^4 -
30\beta^2 +1}{432\beta^4 - 120\beta^2 + 3} &
-\frac{3\beta^2}{36\beta^2 -1} & -\frac{24\beta^3 - 2\beta}{432\beta^4
-120\beta^2 +3}\\
-\frac{3\beta^2}{36\beta^2-1} & -\frac{\beta}{72\beta^2 -2} &
-\frac{3\beta^2}{36\beta^2-1} & \frac{12\beta^2-1}{288\beta^2-8} &
-\frac{\beta}{72\beta^2-2}\\
\frac{16\beta^3}{144\beta^4-40\beta^2+1} &
\frac{8\beta^2}{432\beta^4-120\beta^2+3} & -\frac{24\beta^3
-2\beta}{432\beta^4-120\beta + 3} & -\frac{\beta}{72\beta^2-2} &
-\frac{20\beta^2-1}{864\beta^4 - 240\beta^2+ 6}
\end{pmatrix}
$$

For example, the generating function $$h^{\tau_4}_{(4),(4)}(\beta) =
\frac{6\cdot 6}{4!}\frac{1-20\beta^2}{864\beta^4-240\beta^2+6} = \sum_{k=0}^\infty \big( \frac{36^k}{8} + \frac{4^{k-1}}{2}\big)\beta^{2k}.$$

Let us try to compare this generating function with the known result. Double
Hurwitz numbers
with $\mu= (d)$ are known as one-part double Hurwitz numbers. Consider the
following function
$$\varsigma (z) = e^{z/2} - e^{-z/2} = 2\ \mathrm{sinh}(z/2).$$

A formula for one-part double Hurwitz numbers was derived in.
Considering the particular case $\mu = \nu = (d)$, we use for the comparison a
formula derived in ~\cite{GJV} and reproved in~\cite{J}. In our notation,
this formula reads:
$$ \sum_{r=0}^\infty \frac{z^r}{r!} h^{\tau_4}((d),(d);r) = \frac 1{d^2} \frac
{\varsigma
(d^2z)}{\varsigma (dz)}.$$

Substituting $d=4$, and applying the inverse Borel transform,
$$\frac 1{4^2} \int_0^\infty e^{-t}\frac {\varsigma
(16zt)}{\varsigma (4zt)}\ dt$$
we obtain a perfect agreement with our result.

\end{document}